\documentclass[11pt]{article}
\usepackage{t1enc}
\usepackage[latin1]{inputenc}
\usepackage[english]{babel}
\usepackage{amsmath,amsthm}
\usepackage{amsfonts}
\usepackage{latexsym}
\usepackage[dvips]{graphicx}
\usepackage{graphicx}
\usepackage{color}

 \usepackage{mathrsfs}

\usepackage{amsmath}
\usepackage{amssymb}
\usepackage{amsbsy}

\usepackage{amsfonts}

\usepackage{amsfonts,srcltx,mathrsfs}

\title{Some useful lemmas on  the edge Szeged index}
\author{ Shengjie  He$^1$\footnote{Corresponding author.
Emails: he1046436120@126.com (Shengjie  He)}\\
{\small\em 1. Department of Mathematics, Beijing Jiaotong University, Beijing,
100044, China}\\
  }

\date{} \textwidth 16cm \textheight 22cm \topmargin 0 cm \hoffset
-1.5 cm \voffset 0cm

\newtheorem{theorem}{Theorem}[section]
\newtheorem{lemma}[theorem]{Lemma}

\begin{document}
\baselineskip 0.50cm \maketitle

\begin{abstract}
The edge Szeged index of a graph $G$ is defined as $Sz_{e}(G)=\sum\limits_{uv\in E(G)}m_{u}(uv|G)m_{v}(uv|G)$, where $m_{u}(uv|G)$ (resp., $m_{v}(uv|G)$) is the number of edges whose distance to vertex $u$ (resp., $v$) is smaller than the distance to vertex $v$ (resp., $u$), respectively.
Recently, Wang et al. [G. Wang, S. Li, D. Qi, H. Zhang, On the edge-Szeged index of unicyclic graphs with given diameter, Appl. Math. Comput. 336 (2018) 94-106]
identified the minimum edge-Szeged index of the $n$-vertex unicyclic graphs with a given diameter.
By using another approach, we characterize the graph with minimum edge Szeged index among all the unicyclic graphs with given order and diameter in this paper.

\vspace{0.2cm}

\vspace{0.2cm}

\end{abstract}

\section{Introduction}


Throughout this paper, all graphs we considered are finite, undirected, and simple. Let $G$ be a connected graph with vertex set $V(G)$ and edge set $E(G)$. For a vertex $u \in V(G)$, the degree of $u$, denote by $d_{G}(u)$, is the number of vertices
which are adjacent to $u$. Call a vertex $u$ a pendent vertex of $G$, if $d_{G}(u)=1$ and call an edge $uv$ a pendent edge of $G$, if $d_{G}(u)=1$ or $d_{G}(v)=1$.
For any two vertices $u, v \in V(G)$, let $d(u, v|G)$ denote the distance between $u$ and $v$ in $G$. The diameter of a graph $G$ is the maximum distance between pairs of vertices of $G$. Denote by $P_n$, $S_n$ and $C_n$ a path, star and cycle on $n$ vertices, respectively. Let $\mathbb{N}$ represent the set of natural numbers. If $a$ and $b$ are two natural numbers with $a \leq b$, $[a, b]$ denote the set
$\{  n \in \mathbb{N}|a \leq n \leq b \}$.

The Wiener index is one of the oldest and the most thoroughly studied topological indices. The Wiener index of a graph $G$ is defined as
$$W(G)=\sum\limits_{\{ u, v\} \subseteq V(G) } d(u, v|G).$$
For any edge $e=uv$ of $G$, $V(G)$ can be partitioned into three sets by comparing with the distance of the vertex in $V(G)$ to $u$ and $v$, and the three sets are as follows:
\begin{eqnarray*}N_{u}(e|G)&=&\{ w\in V(G): d(u, w|G) < d(v, w|G) \},\\
N_{v}(e|G)&=&\{ w\in V(G): d(v, w|G) < d(u, w|G) \},\\
N_{0}(e|G)&=&\{ w\in V(G): d(u, w|G) = d(v, w|G) \}.\end{eqnarray*}
The number of vertices of $N_{u}(e|G)$, $N_{v}(e|G)$ and $N_{0}(e|G)$ are denoted by $n_{u}(e|G)$, $n_{v}(e|G)$ and $n_{0}(e|G)$, respectively. If $G$ is a tree, then the formula $W(G)=\sum\limits_{e=uv \in E(G)}n_{u}(e|G)n_{v}(e|G)$
gives a long time known property of the Wiener index. According to this result, a new topological index, named by Szeged index, was introduced by Gutman \cite{Gut.A}, which is an extension of the Wiener index and defined by
$$Sz(G)=\sum\limits_{e=uv \in E(G)}n_{u}(e|G)n_{v}(e|G).$$
For other results on the Wiener and Szeged indices, we refer to \cite{Do.RI,Kh.P,zhang.H,ZhouB.X}.

An edge version of the Wiener index (resp., Szeged index) is also introduced, named edge Wiener index (resp., edge Szeged index). If $e=uv$ and $f=ab$ are two edges of $G$ and $w$ is a vertex of $G$, then the distance between $e$ and $w$ is defined as $d(e,w|G) = {\rm{min}} \{ d(u,w|G),d(v,w|G) \}$ and the distance between $e$ and $f$ is defined as $d(e,f|G) = {\rm{min}} \{ d(u,f|G),d(v,f|G) \}$.
The edge Wiener index \cite{EWiener} is defined as follows:
$$W_{e}(G)= \sum\limits_{ \{e, f \} \subseteq E(G)} d(e, f|G).$$
For $e=uv \in E(G)$, let $M_u(e|G)$ be the set of edges whose distance to the vertex $u$ is smaller than the distance to the vertex $v$, and $M_v(e|G)$ be the set of edges whose distance to the vertex $v$ is smaller than the distance to the vertex $u$. Set $m_u(e|G)=|M_u(e|G)|$ and $m_v(e|G)=|M_v(e|G)|$. Gutman and Ashrafi \cite{Gut.A.R} introduced an edge version of the Szeged index, named edge Szeged index. The edge Szeged index of $G$ is defined as
$$Sz_{e}(G)=\sum\limits_{uv \in E(G)}m_{u}(uv|G)m_{v}(uv|G).$$

In \cite{Gut.A.R}, some basic properties of the edge Szeged index were established by Gutman. It can be checked that the pendent edges make no contributions to the edge Szeged index of a graph. Nadjafi-Arani et al. \cite{EDGE.Relation} proved that for every connected graph $G$, $Sz_{e}(G) \geq W_{e}(G)$ with equality if and only if $G$ is a tree.
 A unicyclic graph of order $n$ is a connected graph with $n$ edges. In \cite{Cai.X}, Cai and Zhou determined the $n$-vertex unicyclic graphs with the largest, the second largest, the smallest and the second smallest edge Szeged indices, respectively. Li \cite{LJP} established a relation between the edge Szeged index and Szeged index of the unicyclic graphs.


Let $\mathcal{U}_{n, d}$ be the set of unicyclic graphs of order $n$ with diameter $d$ for $2 \leq d \leq n-2$.
Tan \cite{Tan} and Shi \cite{Shi} independently determined the graph in $\mathcal{U}_{n, d}$ with minimum Wiener index.
Liu et al. \cite{Yu} characterized the graph in $\mathcal{U}_{n, d}$ with minimum Szeged index.
Recently, Wang et al. \cite{WLI} identified the minimum edge-Szeged index of the $n$-vertex unicyclic graphs with a given diameter.
By using another approach, we determine the graph in $\mathcal{U}_{n, d}$ with  minimum edge Szeged index in this paper. In order to state our results, we need more notations.


Suppose that $l_{1}$, $l_{2}$ and $a$ are three nonnegative integers. Let $P' =P_{l_{1}+1}$ be a path on $l_{1}+1$ vertices with a terminal vertex $u$ and $P'' =P_{l_{2}+1}$ be a path on $l_{2}+1$ vertices with a terminal vertex $v$. Let $S'=S_{a+1}$ be a star on $a+1$ vertices with center vertex $w$. We denote by $P^{a}_{l_{1},l_{2}}$ the tree obtained
from $P'$, $P''$ and $S'$ by identifying $u$, $v$ and $w$ to $u'$, and call $u'$ the root vertex of $P^{a}_{l_{1},l_{2}}$ (see Fig. 1). For
convenience, we write $P^{a}_{l_{1}}$ for $P^{a}_{l_{1},0}$ (see Fig. 1). Note that $P^{a}_{0,0}\cong S_{a+1}$.

\begin{center}   \setlength{\unitlength}{0.7mm}
\begin{picture}(30,40)

\put(-10,25){\circle*{1.5}}
\put(-20,20){\circle*{1.5}}
\put(-20,30){\circle*{1.5}}
\put(-25,17.5){\circle*{1.5}}
\put(-25,32.5){\circle*{1.5}}
\put(-35,12.5){\circle*{1.5}}
\put(-35,37.5){\circle*{1.5}}
\put(-22.5,18.75){\circle*{1}}
\put(-22.5,31.25){\circle*{1}}
\put(-10,25){\line(-2,-1){10}}
\put(-10,25){\line(-2,1){10}}

\put(0,20){\circle*{1.5}}
\put(0,30){\circle*{1.5}}
\put(0,23){\circle*{1}}
\put(0,25){\circle*{1}}
\put(0,27){\circle*{1}}

\put(-10,25){\line(2,-1){10}}
\put(-10,25){\line(2,1){10}}

\put(-25,32.5){\line(-2,1){10}}
\put(-25,17.5){\line(-2,-1){10}}
\put(-50,12.5){\scriptsize$P_{l_{1}+1}$}
\put(-50,37.5){\scriptsize$P_{l_{2}+1}$}

\put(-13,19){\scriptsize$u'$}
\put(1,25){\scriptsize$a$}

\put(60,25){\circle*{1.5}}

\put(50,30){\circle*{1.5}}
\put(45,32.5){\circle*{1.5}}
\put(35,37.5){\circle*{1.5}}
\put(47.5,31.25){\circle*{1}}
\put(60,25){\line(-2,1){10}}

\put(70,20){\circle*{1.5}}
\put(70,30){\circle*{1.5}}
\put(70,23){\circle*{1}}
\put(70,25){\circle*{1}}
\put(70,27){\circle*{1}}

\put(60,25){\line(2,-1){10}}
\put(60,25){\line(2,1){10}}

\put(45,32.5){\line(-2,1){10}}

\put(20,37.5){\scriptsize$P_{l_{1}+1}$}

\put(57,19){\scriptsize$u'$}
\put(71,25){\scriptsize$a$}

\put(55,10){\scriptsize$P_{l_{1}}^{a}$}
\put(-18,10){\scriptsize$P_{l_{1},l_{2}}^{a}$}
\put(-15,0){\scriptsize$\rm{Fig.\,1}$ Graphs $P_{l_{1},l_{2}}^{a}$ and $P_{l_{1}}^{a}$}

\end{picture} \end{center}

Let $C_{3}(P^{n-d-2}_{ \lceil \frac{d-1}{2} \rceil  }, P_{  \lfloor \frac{d-1}{2} \rfloor +1}, S_{1} )$ be the graph obtained from the cycle $C_{3}=v_{1}v_{2}v_{3}v_{1}$
by identifying the root vertex of $P^{n-d-2}_{\lceil  \frac{d-1}{2}\rceil}$ with $v_{1}$, and identifying one end vertex of $P_{\lfloor  \frac{d-1}{2} \rfloor+1}$ with $v_{2}$ (see Fig. 2).

\begin{center}   \setlength{\unitlength}{0.7mm}
\begin{picture}(30,30)

\put(0,25){\circle*{1.5}}
\put(-10,25){\circle*{1.5}}
\put(-20,25){\circle*{1.5}}
\put(-30,25){\circle*{1.5}}
\put(-40,25){\circle*{1.5}}
\put(-50,25){\circle*{1.5}}
\put(-60,25){\circle*{1.5}}
\put(-70,25){\circle*{1.5}}
\put(-50,25){\line(1,0){10}}
\put(-40,25){\line(1,0){10}}
\put(-30,25){\line(1,0){10}}
\put(-10,25){\line(1,0){10}}
\put(-70,25){\line(1,0){10}}

\put(-13,25){\circle*{1}}
\put(-15,25){\circle*{1}}
\put(-17,25){\circle*{1}}
\put(-53,25){\circle*{1}}
\put(-55,25){\circle*{1}}
\put(-57,25){\circle*{1}}
\put(-35,35){\circle*{1.5}}
\put(-30,25){\line(-1,2){5}}
\put(-40,25){\line(1,2){5}}

\put(-35,15){\circle*{1.5}}
\put(-25,15){\circle*{1.5}}
\put(-30,25){\line(-1,-2){5}}
\put(-30,25){\line(1,-2){5}}
\put(-28,15){\circle*{1}}
\put(-30,15){\circle*{1}}
\put(-33,15){\circle*{1}}

\put(-74,27){\scriptsize$u_{0}$}
\put(-64,27){\scriptsize$u_{1}$}
\put(-14,27){\scriptsize$w_{1}$}
\put(-4,27){\scriptsize$w_{0}$}
\put(-50,30){\scriptsize$u_{\lfloor \frac{d-1}{2} \rfloor}$}
\put(-35,30){\scriptsize$w_{\lceil \frac{d-1}{2} \rceil}$}
\put(-23,15){\scriptsize$n-d-2$}

\put(-65,7){\scriptsize$C_{3}(P^{n-d-2}_{ \lceil \frac{d-1}{2} \rceil  }, P_{  \lfloor \frac{d-1}{2} \rfloor +1}, S_{1} )$}

\put(90,25){\circle*{1.5}}
\put(80,25){\circle*{1.5}}
\put(70,25){\circle*{1.5}}
\put(60,25){\circle*{1.5}}
\put(50,25){\circle*{1.5}}
\put(40,25){\circle*{1.5}}
\put(30,25){\circle*{1.5}}
\put(20,25){\circle*{1.5}}
\put(40,25){\line(1,0){10}}
\put(50,25){\line(1,0){10}}
\put(60,25){\line(1,0){10}}
\put(80,25){\line(1,0){10}}
\put(20,25){\line(1,0){10}}

\put(77,25){\circle*{1}}
\put(75,25){\circle*{1}}
\put(73,25){\circle*{1}}
\put(37,25){\circle*{1}}
\put(35,25){\circle*{1}}
\put(33,25){\circle*{1}}

\put(55,15){\circle*{1.5}}
\put(65,15){\circle*{1.5}}
\put(60,25){\line(-1,-2){5}}
\put(60,25){\line(1,-2){5}}
\put(62,15){\circle*{1}}
\put(60,15){\circle*{1}}
\put(57,15){\circle*{1}}

\put(16,27){\scriptsize$u_{0}$}
\put(26,27){\scriptsize$u_{1}$}
\put(73,27){\scriptsize$u_{d-1}$}
\put(86,27){\scriptsize$u_{d}$}

\put(67,15){\scriptsize$n-d-1$}

\put(56,29){\scriptsize$u_{\lfloor \frac{d}{2} \rfloor}$}

\put(55,7){\scriptsize$T_{n, d, \lfloor \frac{d}{2} \rfloor}$}

\put(-35,-2){\scriptsize$\rm{Fig.\,2}$ Graphs $C_{3}(P^{n-d-2}_{ \lceil \frac{d-1}{2} \rceil  }, P_{  \lfloor \frac{d-1}{2} \rfloor +1}, S_{1} )$ and $T_{n, d, \lfloor \frac{d}{2} \rfloor}$}

\end{picture} \end{center}

Let $G$ be a unicyclic graph of order $n$ with diameter $d$. Then $1 \leq d \leq n-2$. If $d=1$, then $G\cong C_{3}$. If $d=2$, then $G\in \{ C_{4}, C_{3}(P_{1}^{0}, S_{1}, S_{1})\}$ for $n=4$; $G\in \{ C_{5}, C_{3}(P_{1}^{1}, S_{1}, S_{1})\}$ for $n=5$ and $G= C_{3}(P_{1}^{n-4}, S_{1}, S_{1})$ for $n \geq 6$.
For $1 \leq d \leq 2$, it is not difficult to determine $Sz_e(G)$. Hence throughout the paper, we assume that $n \geq 6$ and $3 \leq d \leq n-2$. Our main result is the following theorem.

\begin{theorem}\label{Th1} Among the unicyclic graphs of order $n \geq 6$ with diameter $d$ for $3 \leq d \leq n-2$, $C_{3}(P^{n-d-2}_{ \lceil \frac{d-1}{2} \rceil  }, P_{  \lfloor \frac{d-1}{2} \rfloor +1}, S_{1} )$ is the unique graph with minimum edge Szeged index.
\end{theorem}

The rest of this paper is organized as follows. In Section 2, we present some useful lemmas. In Section 3, we establish some transformations of graphs which keep
the diameter and the order of the graphs, but decrease the edge Szeged index of the graphs. In section 4, Theorem \ref{Th1} is proved.

\section{Useful lemmas}


We first introduce a result useful for  computing the edge Szeged index of a
unicyclic graph.
Let $g \geq 3$ be an integer, and let $C_{g}=v_{1}v_{2} \cdots v_{g}v_{1}$ be a cycle on $g$ vertices. Let $T_{1}, T_{2}, \cdots, T_{g}$ be vertex-disjoint trees,
and let the root vertex of $T_{i}$ be $u_{i}$ for $1 \leq i \leq g$. Denote by $C_{T_{1}, T_{2}, \cdots, T_{g}}^{u_{1}, u_{2}, \cdots, u_{g}}$ the unicyclic graph obtained from $C_{g}$ by identifying the root vertex $u_{i}$ of $T_{i}$ with $v_{i}$ for each $1 \leq i \leq g$.
If there is no confusion, we write $C_{g}(T_{1}, T_{2}, \cdots, T_{g})$ for $C_{T_{1}, T_{2}, \cdots, T_{g}}^{u_{1}, u_{2}, \cdots, u_{g}}$. Any $n$-vertex  unicyclic graph $G$ with a  $g$-cycle is of the form $C_{g}(T_{1}, T_{2}, \cdots, T_{g})$, where $\sum_{i=1}^{g}|T_{i}|=n=\sum_{i=1}^{g}|E(T_{i})|+g$.
For $v \in V(G)$, we define $$D(v|G)=\sum_{w\in V(G)}d(v, w|G).$$ For an integer $g$, define $$ \delta(g)=\left\{
                                                 \begin{array}{ll}
                                                   1, & \hbox{if $g$ is odd ;} \\
                                                   0, & \hbox{if $g$ is even.}
                                                 \end{array}
                                               \right.
  $$

\begin{lemma}\label{lem2.8}{\rm \cite{LJP}} Let $G=C_{r}(T_{1}, T_{2}, \cdots, T_{r})$ with $|V(G)|=n$. Then

$$ Sz_{e}(G)= Sz(G)+ \sum\limits_{i=1}^{r}D(v_{i}|T_{i})-n^{2}+ \left\{
                                                 \begin{array}{ll}
                                                   nr & \hbox{if $r$ is odd ;} \\
                                                   r & \hbox{if $r$ is even.}
                                                 \end{array}
                                               \right.
  $$

\end{lemma}




The following are some lemmas which will be used in our proof.

\begin{lemma}\label{lem2.3}{\rm \cite{S.H.H.Y}}
Let $G$ and $G'$ be the graphs shown as in $Fig. \, 3$, where $G$ consists of $G_0$ and $G_1$ with a common vertex $u$, and $G'$ consists of $G_0$ and $G_2$ with a common vertex $u$. Then each of the followings holds:\\
{\em(i)} For any edge $e=w_1w_2\in E(G_0)$ and $1\leq i\leq 2$, we have
                 $$m_{w_i}(e|G)=m_{w_i}(e|G_0)+\tau(u)|E(G_1)|,$$
                 where $$ \tau(u)=\left\{
                                                 \begin{array}{ll}
                                                   1, & \hbox{$u\in N_{w_i}(e|G_0)$ ;} \\
                                                   0, & \hbox{otherwise.}
                                                 \end{array}
                                               \right.
  $$
{\em(ii)} If $|E(G_1)|=|E(G_2)|$, then
\begin{eqnarray*}
\sum\limits_{e=w_1w_2 \in E(G_0)} m_{w_1}(e|G)m_{w_2}(e|G)&=&\sum\limits_{e=w_1w_2 \in E(G_0)} m_{w_1}(e|G')m_{w_2}(e|G').
\end{eqnarray*}
\end{lemma}

\begin{center}   \setlength{\unitlength}{0.7mm}
\begin{picture}(30,40)

\put(-15,30){\circle*{1}}
\put(-25,30){\circle{20}}
\put(-5,30){\circle{20}}
\put(-18,5){\scriptsize$G$}
\put(42,5){\scriptsize$G'$}

\put(-21,30){\scriptsize$u$}
\put(39,30){\scriptsize$u$}

\put(45,30){\circle*{1}}

\put(-27,15){\scriptsize$G_{0}$}
\put(-7,15){\scriptsize$G_{1}$}
\put(33,15){\scriptsize$G_{0}$}
\put(53,15){\scriptsize$G_{2}$}

\put(35,30){\circle{20}}
\put(55,30){\circle{20}}

\put(-15,-2){\scriptsize Fig. 3. $G$ and $G'$ in Lemma 2.2}
\end{picture} \end{center}

\begin{lemma}\label{lem2.5}{\rm \cite{S.H.H.Y}} Let $G$ be a graph of order $n$ with a cycle $C_g=v_1v_2 \cdots v_{g}v_1$.
Assume that $G-E(C_{g})$ has exactly $g$ components $G_1,G_2,\cdots, G_g$, where $G_i$ is the component of $G-E(C_{g})$ that contains $v_i$ for $1\leq i\leq g$.
Let
 \begin{eqnarray*}G'=G- \displaystyle\cup_{i=2}^g\{wv_i: w\in N_{G_i}(v_i)\}+ \cup_{i=2}^g\{wv_1:  w\in N_{G_i}(v_i)\}.\end{eqnarray*}
Then $Sz_e(G') \leq Sz_e(G)$ with equality if and only if $C_{g}$ is an end-block, that is, $G\cong G'$.
\end{lemma}

\begin{lemma}\label{lem2.2}{\rm \cite{Gut.A.R}}  If $T$ is a tree with order $n$, then
$Sz_{e}(T)=Sz(T)-(n-1)^{2}.$
\end{lemma}

Note that if $T$ is a tree, $W_{e}(T)=Sz_e(T)$ and $W(T)=Sz(T)$. The conclusion of  Lemma \ref{lem2.2} can also be written as $W_{e}(T)=W(T)-(n-1)^{2}.$

Let $n\geq 3$ and $2 \leq d \leq n-1$ be integers. Let $P_{d+1}=u_{0}u_{1} \cdots u_{d}$ be a path of order $d+1$ and $T_{n, d, \lfloor \frac{d}{2} \rfloor }$
be the tree obtained from $P_{d+1}$ by attaching $n-d-1$ pendent vertices to $u_{\lfloor \frac{d}{2} \rfloor }$ (see Fig. 2).
Note that $T_{n, d, \lfloor \frac{d}{2} \rfloor } \cong P^{n-d-1}_{\lfloor \frac{d}{2} \rfloor , \lceil \frac{d}{2} \rceil }$.

\begin{lemma}\label{lem2.6}{\rm \cite{HQLIU}}
Among the trees of order $n$ with diameter $2 \leq d\leq n-2$, $T_{n, d, \lfloor \frac{d}{2} \rfloor }$ has minimum Wiener index.
\end{lemma}

\begin{lemma}\label{lem2.9}{\rm \cite{Yu}} Let $C_g=v_1v_2 \cdots v_{g}v_1$ and $d_{i, j}=d(v_{i}, v_{j}|C_g)$ for $i, j \in \{1, 2, \cdots, g \}$.
\\
{\em(i)} If $g$ is an even number,
$$  (d_{2, j}-d_{1, j}+1, d_{g, j}-d_{1, j}+1) =     \left\{
  \begin{array}{ll}
    (0, 2), & \hbox{if $2 \leq j \leq \frac{g}{2}$;} \\
    (0, 0), & \hbox{if $j=\frac{g}{2}+1$;} \\
    (2, 0), & \hbox{if $\frac{g}{2}+2 \leq j \leq g-1.$}
  \end{array}
\right.
$$
{\em(ii)} If $g$ is an odd number,
$$      (d_{2, j}-d_{1, j}+1, d_{g, j}-d_{1, j}+1) =\left\{
  \begin{array}{ll}
    (0, 2), & \hbox{if $2 \leq j \leq \frac{g-1}{2}$;} \\
    (0, 1), & \hbox{if $j=\frac{g+1}{2}$;} \\
    (1, 0), & \hbox{if $j=\frac{g+3}{2}$;} \\
    (2, 0), & \hbox{if $\frac{g+5}{2} \leq j \leq g-1.$}
  \end{array}
\right.
$$
\end{lemma}

\section*{Acknowledgments}

This work was supported by the National Natural Science Foundation of China
(Nos. 11371052,11\\771039), the Fundamental Research Funds for the Central Universities (Nos. 2016JBM071, 2016JBZ012) and the $111$ Project of China (B16002).


\end{document}